\documentclass[11pt]{article}
\usepackage{amsfonts}
\usepackage{color,xcolor}

\begin{document}

\title{ Analysis of a mathematical model for tumor growth with 
  Gibbs-Thomson relation}

\author{Junde Wu}
\date{{\small Department of Mathematics, Soochow University,
  Suzhou, \\
  Jiangsu 215006, PR China,  E-mail: wujund@suda.edu.cn}\\ [0.1cm]
  }
\maketitle

\begin{abstract}
  In this paper we study a mathematical model for the growth of  
  nonnecrotic solid tumor. The tumor is assumed to be radially 
  symmetric and its radius $R(t)$ is an unknown function of time 
  $t$ as tumor growth, and the model is in the form of a free boundary
  problem. The feature of the model is that a Gibbs-Thomson relation 
  is taken into account, which resulting an interesting phenomenon 
  that there exist two stationary solutions (depending on the model 
  parameters). The global existence and uniqueness of solution are 
  established. By denoting $c$ the ratio of the diffusion time scale 
  to the tumor doubling time scale, we prove that for sufficiently 
  small $c>0$, the stationary solution with the larger radius is 
  asymptotically stable, and the other smaller one is unstable.
\medskip

\medskip

\medskip

\end{abstract}

\section{Introduction}
\setcounter{equation}{0}\hskip 1em

  In the last several decades, great attention has been attracted to 
  mathematical models of tumor growth for their own both biological 
  and mathematical interests, cf.  [\ref{fri-07}, \ref{gre-72}--\ref{roose-07}]
  and references therein. Mathematical analysis of these models 
  can help us understanding the mechanism of tumor growth and 
  accessing tumor treatment strategy.  On the other hand, a lot of
  mathematical challenges arise in tumor models, and many interesting 
  and illuminative results have been established, cf. [\ref{cui-02}--\ref{fri-rei-99},
  \ref{wu-cui-09}, \ref{wu-zhou-13}, \ref{zhou-wu-15}] and references therein. 
  
  In this paper we study a tumor model in the form of a free boundary problem.
  Since solid tumors grow with spheroid-shaped,  tumor region is assumed to
  be a spheroid with radius $R(t)$ at time $t>0$, the proliferation of tumor cells  
  is assumed to be dependent only on located concentration of nutrient $\sigma(r,t)$, 
  which is diffusing 
  within tumor region, and tumor growth is governed by the mass 
  conservation law. The model considered here is given as follows:
\begin{equation}\label{1.1}
   c{\partial \sigma\over \partial t}= {1\over r^2}{\partial\over\partial r}
   \Big(r^2{\partial\sigma\over\partial r}\Big)-\lambda\sigma  \qquad 
    \mbox{for}\;\;r<R(t),\quad t>0,
\end{equation}
\begin{equation}\label{1.2}
  {\partial \sigma \over \partial r}(0,t)=0,\qquad
  \sigma(R(t),t)=G(t)
   \qquad  \mbox{for}\;\; t>0,
\end{equation}
\begin{equation}\label{1.3}
   {d R\over dt}=
   {1\over R^2}\int_0^{R} \mu(\sigma-\tilde\sigma) r^2 dr \qquad \mbox{for}\;\; t>0, 
\end{equation}
\begin{equation}\label{1.4}
   \sigma(r,0)=\sigma_0(r)  \qquad\mbox{for}\;\; 0\le r\le R_0,
\end{equation}
\begin{equation}\label{1.5}
      R(0)=R_0.
\end{equation}
 where $c$, $\lambda$, $\tilde\sigma$, $\mu$ are positive
 dimensionless constants, among which $c$ is the ratio of the 
 nutrient diffusion time scale  ($\sim$ minutes) to the tumor 
 doubling time scale ($\sim$ days), and $c\ll 1$; $\lambda$ 
 is the nutrient consumption rate;  $\tilde\sigma$ is a threshold 
 value of nutrient concentration for apoptosis;  $\mu$ is the 
 proliferation rate of tumor cells;  $G(t)$ is a given function 
 representing the external nutrient supply; $\sigma_0(r)$ 
 and $R_0$ are the given initial data.
 
 In Byrne and Chaplain [\ref{byr-cha-95}], the external 
 nutrient concentration is assumed to be constant $\bar\sigma$,
 i.e., $G(t)=\bar\sigma$.  In this case, Friedman and Reitich
 [\ref{fri-rei-99}] proved that the tumor model  (\ref{1.1})--(\ref{1.5}) 
 has a unique radially symmetric stationary solution for
 $0<\tilde\sigma<\bar\sigma$, and it is asymptotically stable 
 for sufficiently small $c$. Later, Cui [\ref{cui-02}] 
 extended the result to the inhibitor-presence case and Cui [\ref{cui-05}] 
 further investigated the above model with general nutrient consumption 
 function and cell proliferation function. Recently, Xu [\ref{xu-16}] 
 considered the case that $G(t)$ is given by a periodic function, 
 global well-posedness  and some asymptotic behavior of solutions 
 were derived.
 
  One disadvantage of above assumptions on $G(t)$ is that though the nutrient 
  is continuous across the tumor boundary, but the flux of nutrient is not.  
  By contrast, Byrne and Chaplain [\ref{byr-cha-96}] assumed that 
  energy is expended in maintaining the tumor's compactness by cell-to-cell 
  adhesion on tumor boundary and the nutrient acts as a source of energy, 
  so the nutrient concentration on the tumor boundary is less than 
  the external supply $\bar\sigma$, and the difference satisfies a 
  Gibbs-Thomson relation, i.e., the difference of nutrient 
  concentration across the tumor boundary $r=R(t)$ is proportional to the 
  mean curvature which is given by $1/R(t)$. More precisely, Byrne and Chaplain 
  [\ref{byr-cha-96}] assumed that $G(t)=\bar\sigma(1-\gamma/R(t))$, where 
  $\gamma$ is a positive constant representing the cell-to-cell adhesiveness.
  In quasi-stationary case $c=0$ and replacing $\lambda\sigma$ by $\lambda$ 
  of equation (\ref{1.1}), Byrne and Chaplain 
   [\ref{byr-cha-96}] studied existence and uniqueness of solution and 
   the linear stability of stationary solutions, numerical verification was 
   also performed. 
   
   Roose, Chapman and Maini [\ref{roose-07}] pointed out that the tumor model
   (\ref{1.1})--(\ref{1.5}) with $G(t)=\bar\sigma(1-\gamma/R(t))$, which is induced 
   by Gibbs-Thomson relation,  has a number of interesting points. Though it seems
   speculative, but it may be possible to check its veracity in experiment. It is 
   significant to analyze how the Gibbs-Thomson relation effects the growth 
   of tumors, which can be also tested in experiments and clinical laboratory. 
   Note that for any positive constant $\gamma$, if $R(t)<\gamma$,
   then $G(t)= \bar\sigma(1-\gamma/R(t))<0$. It is unreasonable since the nutrient 
   concentration must be always nonnegative. For this reason, we introduce a 
   simple modification and let 
\begin{equation}\label{1.6}
  G(t)=\bar\sigma(1-\gamma/R(t))H(R(t)),
\end{equation}  
  where $H(\cdot)$ is a smooth function such that 
  $H(r)=0$ for $r\le \gamma$, $H(r)=1$ for $r\ge 2\gamma$, and $0\le H'(r)\le 1/\gamma$. 
  In this paper, we shall make a rigorous analysis of problem (\ref{1.1})--(\ref{1.5}) with
  $G(t)$ given by (\ref{1.6}), and study the effect of 
  Gibbs-Thomson relation.
  
  In Section 2, we shall prove the global existence and uniqueness of solution,
  based on a priori estimate and fixed point method. In Section 3, we study the
  quasi-stationary case $c=0$. We shall prove that there may exist two, or a 
  unique, or none radially symmetric stationary solutions depending on model 
  parameters.  It is interesting that there may exist two radially symmetric 
  stationary solutions, 
  which is different from the uniqueness of stationary solution for constant 
  $G(t)=\bar\sigma$ in [\ref{fri-rei-99}] and periodic function $G(t)$ in 
  [\ref{xu-16}]. By the linearized stability principle, 
  we shall see that in quasi-stationary case $c=0$, the radially symmetric
  stationary solution with the larger radius is asymptotically stable and the 
  other smaller one is unstable.
  
  In Section 4 we study fully non-stationary case $c>0$. By using a comparison
  method and some analysis techniques motivated by [\ref{cui-02}, \ref{cui-05}, 
  \ref{fri-rei-99}], we shall prove that for sufficiently small $c>0$, the same stability 
  results still hold as the quasi-stationary case
  $c=0$.  In the last section, we make a conclusion and give some interesting 
  biological implications.

\medskip 
\hskip 1em

\section{Global existence and uniqueness}
\setcounter{equation}{0} 
\hskip 1em

  In this section we study global existence and uniqueness of problem (\ref{1.1})--(\ref{1.6}).
  Throughout this paper we assume that the initial data $\sigma_0(r)$ and $R_0$ satisfy 
  the following conditions:
  \begin{equation}\label{2.1}
  R_0>0;\quad \sigma_0(r)\in C^2[0,R_0], \quad 0\le \sigma_0(r) \le \bar\sigma,
  \quad
  \sigma_0'(0)=0\;\; \mbox{ and }\;\; \sigma_0(R_0)=G(0).
\end{equation}

\medskip

  {\bf  Theorem 2.1} \ \  Problem (\ref{1.1})--(\ref{1.6}) has a 
  unique solution $(\sigma(r,t), R(t))$ for all $t>0$,
  and there hold following assertions:
\begin{equation}\label{2.2}  
  0\le \sigma(r,t)\le \bar\sigma\qquad \mbox{for}\;\; 0\le r\le R(t),\;\; t\ge 0,
\end{equation}
\begin{equation}\label{2.3}
 -{1\over3} \mu\tilde\sigma \le{R'(t)\over R(t)} \le 
 {1\over3}\mu(\bar\sigma-\tilde\sigma)\qquad \mbox{for}\;\; t\ge0,
\end{equation}
\begin{equation}\label{2.4}
 R_0\exp(-{1\over3}\mu\tilde\sigma t)<R(t)\le R_0 \exp({1\over3}\mu(\bar\sigma-\tilde\sigma) t) \qquad\mbox{for}\;\; 
 t\ge0.
\end{equation}
\medskip

 {\bf Proof}. \ \  We first assume that $(\sigma(r,t), R(t))$ is a solution of 
 problem (\ref{1.1})--(\ref{1.6}). By the maximum principle, we immediately 
 have $0\le\sigma(r,t)\le\bar\sigma$ for $0\le r\le R(t)$, $t>0$. 
 
   By (\ref{1.3}),
$$
  {d R\over d t}={1\over R^2(t)}\int^{R(t)}_0 \mu (\sigma(r,t)-\tilde\sigma) r^2 dr,
$$
  we have
$$
 -{1\over 3}\mu\tilde\sigma \le {R'(t)\over R(t)}\le {1\over 3}\mu(\bar\sigma-\tilde\sigma),
$$
  then (\ref{2.3}) and (\ref{2.4}) follow obviously. 
  
  Next, we prove the existence and uniqueness of solution to the problem. 
  For arbitrary $T>0$, we introduce a metric space $(M_T,d)$ as follows:
  The set $M_T$ consists of vector functions $(\sigma(r,t),R(t))$ which
  satisfy 
  
  $(i)$ $R\in C[0,T]$, $R(0)=R_0$, and 
$$
  R_0\exp(-{1\over3}\mu\tilde\sigma t)\le R(t)\le R_0\exp({1\over 3} \mu(\bar\sigma-\tilde\sigma)t)\qquad\mbox{for}\;\;
  0<t\le T.
$$
  
  $(ii)$ $\sigma\in C([0,\infty)\times [0,T])$, and
$$
\begin{array}{l}
 \displaystyle 0
  \le \sigma(r,t)\le \bar\sigma\qquad \mbox{for}\;\; 0\le r\le R(t),\;\;
  0<t\le T,
  \\ [0.3 cm]
 \displaystyle \sigma(R(t),t)=\bar\sigma(1-{\gamma\over R(t)})H(R(t)) 
 \qquad\mbox{for}\;\;0<t\le T,
  \\ [0.3 cm]
  \sigma(r,0)=\sigma_0(r)\qquad\mbox{for}\;\; 0\le r\le R_0.
\end{array}
$$
  The metric $d$ is defined by
$$
  d((\sigma_1,R_1),(\sigma_2,R_2))=\max_{r\ge0,0\le t\le T} 
  |\sigma_1(r,t)-\sigma_2(r,t)|
  +\max_{0\le t\le T} |R_1(t)-R_2(t)|.
$$
  It is clear that $(M_T, d)$ is a complete metric space. For a given $(\sigma,R)\in M_T$, 
  let $\hat R$ be the solution of the following initial value problem
$$
\left\{
\begin{array}{l}
  \displaystyle {d \hat R\over d t}={\hat R\over R^3}\int_0^R \mu(\sigma-\tilde\sigma)r^2 dr
  \qquad\mbox{for}\;\; 0\le t\le T,
  \\ [0.5 cm]
  \hat R(0)=R_0.
\end{array}
\right.
$$
  Clearly,
$$
  \displaystyle \hat R(t)=R_0\exp(\int_0^t K(\tau)d\tau),\qquad
  K(t)={1\over R^3(t)}\int_0^{R(t)} \mu(\sigma-\tilde\sigma)r^2 dr.
$$
  Since $0\le\sigma(r,t)\le\bar\sigma$, we see that 
  $-{1\over3}\mu\tilde\sigma\le K(t)\le{1\over 3}\mu(\bar\sigma-\tilde\sigma)$
  and so that $\hat R$ satisfies condition $(i)$.  Next, we 
  consider the following problem
\begin{equation}\label{2.5}
\left\{
\begin{array}{l}
   \displaystyle c{\partial \hat \sigma\over \partial t}= {1\over r^2}{\partial\over\partial r}
   \Big(r^2{\partial\hat\sigma\over\partial r}\Big)-\lambda\hat\sigma  \qquad 
    \mbox{for}\;\;0<r<\hat R(t),\;\; 0<t\le T, 
    \\ [0.5 cm]
   \displaystyle {\partial \hat \sigma \over \partial r}(0,t)=0,\qquad 
   \hat\sigma(\hat R(t),t)=\hat G(t)
   \qquad\mbox{for}\;\; 0<t\le T,
    \\ [0.5 cm]
   \hat\sigma(r,0)=\sigma_0(r)\qquad\mbox{for}\;\; 0\le r\le R_0.
\end{array}
\right.
\end{equation}
  where $\hat G(t)=\bar\sigma(1-\gamma/\hat R(t))H(\hat R(t))
  \in C[0,T]$.
  By letting $u(y,t)=\hat\sigma({\hat R(t)y},t)$, 
  it is equivalent to the following problem
\begin{equation}\label{2.6}
\left\{
\begin{array}{l}
   \displaystyle c{\partial u\over \partial t}= {1\over \hat R^2(t)y^2}{\partial\over\partial y}
   \Big(y^2{\partial u\over\partial y}\Big)+{c \hat R'(t)\over \hat R(t)}y{\partial u\over \partial y}
   -\lambda u
   \qquad  \mbox{for}\;\;0<y<1,\; 0<t\le T, 
    \\ [0.5 cm]
   \displaystyle {\partial u \over \partial y}(0,t)=0,\qquad 
   u(1,t)=\hat G(t)
   \qquad\mbox{for}\;\; 0<t\le T,
    \\ [0.5 cm]
   u(y,0)=\sigma_0( R_0 y)
   \qquad\mbox{for}\;\; 0\le y\le 1.
\end{array}
\right.
\end{equation}
  Since all coefficients of the above differential equations are bounded, thus by 
  a standard theory of parabolic equations, we see that there exists a unique 
  solution $u(y,t)\in C([0,1]\times [0,T])$. 
  We extend $u(y,t)$ such that $u(y,t)=\hat G(t)$ for $y\ge 1$, $0\le t\le T$,   
  and get a corresponding $\hat \sigma$. 
  By comparison, we have $0\le \hat \sigma(r,t)\le \bar\sigma$ 
  and condition $(ii)$ is satisfied. Hence for small $T>0$,
  we can define a mapping $W: M_T\to M_T$ such that 
$$
  W(\sigma,R)=(\hat \sigma, \hat R).
$$
  By a similar analysis of [\ref{cui-05}], we can further show that $W$ is a 
  contraction mapping on $M_T$ for sufficiently small $T$. Then by Banach 
  fixed point theorem, we get the local existence and uniqueness of 
  problem (\ref{1.1})--(\ref{1.6}). 
  
  Finally, since (\ref{2.2}) and (\ref{2.3}) do not depend on initial data $\sigma_0(r)$ and 
  $R_0$, we can extend the solution to all $t>0$. \qquad $\Box$

\medskip 
\hskip 1em

\section{Quasi-stationary case $c=0$}
\setcounter{equation}{0} 
\hskip 1em

  In this section we study quasi-stationary case $c=0$ of 
  free boundary problem (\ref{1.1})--(\ref{1.6}).  For simplification 
  of notations, by a rescaling argument, we always set 
  $\lambda=\bar\sigma=1$ later on.  
  
  First, we study radially symmetric stationary solution which is
  denoted by
  $(\sigma_s(r), R_s)$ for $R_s>0$.
  It is easy to see that 
\begin{equation}\label{3.1}
  \sigma_s(r)=(1-{\gamma\over R_s}){R_s\sinh r\over r\sinh R_s}H(R_s).
\end{equation}
  Substituting it into the right term of equation $(\ref{1.3})$ and by using 
  the relation $\displaystyle {d R_s\over dt}=0$, we see $R_s>0$ satisfies
\begin{equation}\label{3.2}
   \Big(1-{\gamma\over R_s}\Big){R_s\coth R_s-1\over R_s^2}H(R_s)
         - {1\over 3}\tilde\sigma=0.
\end{equation}
 Denote 
\begin{equation}\label{3.3}
  f(r):=(1-{\gamma\over r}){r\coth r-1\over r^2}
  \quad\mbox{and}\quad
  F(r):=3f(r)H(r)\quad\mbox{for}\;\;r>0.
\end{equation}
 Thus $R_s>0$ is the root of the equation 
 $ F(r)={\tilde\sigma}$.
 
  It is easy to verify that 
\begin{equation}\label{3.4}  
    F(r)\left\{
  \begin{array}{l}
  =0,\qquad\mbox{for}\;\; 0<r\le\gamma,
  \\ [0.2 cm]
  >0,\qquad\mbox{for}\;\; r>\gamma,
  \end{array}
  \right.
  \qquad\mbox{and}\;\;
  \lim_{r\to+\infty} F(r)=0.
\end{equation}
   Moreover, by the proof of Theorem 1.1 in [\ref{wu-16}],  we see
   that $F(r)$ has a unique extremum point $r_\#\in (2\gamma, 2\gamma+2)$ 
   such that 
\begin{equation}\label{3.5}
  F'(r)\left\{
  \begin{array}{l}
  >0,\qquad\mbox{for}\;\; 0<r<r_\#,
  \\ [0.2 cm]
  =0,\qquad\mbox{for}\;\; r=r_\#,
  \\ [0.2 cm]
  <0,\qquad\mbox{for}\;\; r>r_\#,
  \end{array}
  \right.
\end{equation}
  and
\begin{equation}\label{3.6}
 0<\theta_*:= F(r_\#)=\max_{r>0} F(r)<1.
\end{equation}
   It immediately follows that:
     
  $(i)$ If $\tilde\sigma>\theta_*$, then equation $F(r)=\tilde\sigma$
    has no positive solution;
    
   $(ii)$ If $\tilde\sigma=\theta_*$, then equation
   $F(r)=\tilde\sigma$ has a unique positive solution $R_s=r_\#$; 
   
  $(iii)$ If $0<\tilde\sigma<\theta_*$, then equation
   $F(r)=\tilde\sigma$ has two positive solutions $R_{s1}$ and $R_{s2}$ 
   satisfying $\gamma<R_{s1}<r_\#<R_{s2}$ with $F'(R_{s1})>0$ and
   $F'(R_{s2})<0$.
   
   Note that in case $G(t)\equiv 1$, for $0<\tilde\sigma<1$, there exists 
   a unique radially symmetric stationary solution. While in case 
   $G(t)$ given by (\ref{1.6}), we see that for $0<\tilde\sigma<\theta_*$,
   the problem has two such stationary solutions.
   Later on,
   we focus on this interesting case $0<\tilde\sigma<\theta_*$.
   
   Let $c=0$, for any given function $R(t)\in C^1[0,\infty)$,
   we solve problem (\ref{1.1})--(\ref{1.2}) and get 
\begin{equation}\label{3.7}
  \sigma(r,t)=
   (1-{\gamma\over R(t)}){R(t)\sinh r\over r\sinh R(t)}H(R(t)).
\end{equation}
  By substituting it into (\ref{1.3}) we reduce the free boundary 
  problem into the following equation
\begin{equation}\label{3.8}
   {d R\over dt}={1\over 3}\mu [F(R)-\tilde \sigma]R.
\end{equation}
  Clearly, by classical linearized stability principle of differential equations,
  we can get the stability of radially stationary solutions. In conclusion,
  we have 
  
\medskip

  {\bf Theorem 3.1} \ \ {\em  Let $0<\tilde\sigma<\theta_*$.
  Free boundary problem $(\ref{1.1})$--$(\ref{1.6})$ has two radially symmetric
  stationary solutions with radius $R_{s1}$ and $R_{s2}$, ($R_{s1}<R_{s2}$), 
  respectively. In quasi-stationary case $c=0$, the stationary solution with the larger radius
  $R_{s2}$ is asymptotically stable and the other smaller one with radius $R_{s1}$ is 
  unstable. More precisely, we have
$$
 \lim_{t\to\infty } R(t)=
  \left\{
  \begin{array}{ll}
  0,\quad& \mbox{for}\;\; 0<R_0<R_{s1},
  \\ [0.2 cm]
  R_{s2}, \qquad& \mbox{for}\;\; R_0>R_{s1}.
  \end{array}
  \right. 
$$
 }
  
\medskip

 {\bf Remark 3.2} \ \  We regard $\gamma$ as a variable and discuss the effect
  of Gibbs-Thomson relation on tumor growth. 
  Rewrite $F(r)$, $\theta_*$ and $R_s$ as $F(r,\gamma)$, $\theta_*(\gamma)$ 
  and $R_s(\gamma)$, respectively,  by regarding them as functions 
  depending on $\gamma$. We have
$$
  {\partial F\over \partial \gamma}=-3{r\coth r-1\over r^3}H(r)<0,
  \qquad\mbox{for}\;\; r>\gamma.
$$
  It implies that $\theta_*'(\gamma)<0$ and for $0<\tilde\sigma<\theta_*(\gamma)$,
  there hold $R_{s1}'(\gamma)>0$ and $R_{s2}'(\gamma)<0$. By Theorem 3.1, we see that
  the radius of the stable radially symmetric stationary solution is decreasing on $\gamma$.
  It implies that increasing cell-to-cell adhesiveness may play a positive role
  on making tumor more stable.
  
\medskip

\hskip 1em

\section{Asymptotic behavior and stability}
\setcounter{equation}{0} 
\hskip 1em

  In this section we study asymptotic behavior of solution $(\sigma(r,t),R(t)) $
  to free boundary problem (\ref{1.1})--(\ref{1.6}).

  First, if the concentration of external nutrient supply is less than
  the threshold value for apoptosis, the tumor will starve and shrink to zero. 
  More precisely,  we have
  
\medskip

  {\bf Theorem 4.1} \ \ {\em If $\tilde\sigma>\bar\sigma$, then for any $c>0$ and
  any given initial data $(\sigma_0(r), R_0)$ satisfying $(\ref{2.1})$, we have
$$
  \lim_{t\to\infty}R(t)=0.
$$ 
}

 {\bf Proof}. \ \  By Theorem 2.1 we see that there exists a unique global solution
 $(\sigma(r,t),R(t))$ of free boundary problem (\ref{1.1})--(\ref{1.6}). By (\ref{2.4})
 we have $R(t)\le R_0\exp({1\over3}(\bar\sigma-\tilde\sigma)t)$, since 
 $\bar\sigma<\tilde\sigma$, we obtain that $\lim_{t\to\infty}R(t)=0$.
 \qquad$\Box$
 
\medskip

  Next we focus on the case $0<\tilde\sigma<\theta_*$ where there exist two 
  radially symmetric stationary solutions denoted by $(\sigma_{s1}(r), R_{s1})$
  and $(\sigma_{s2}(r),R_{s2})$, with $R_{s1}<R_{s2}$, respectively. 
  
  Let $(\sigma(r,t),R(t))$ is a solution of problem  (\ref{1.1})--(\ref{1.6}) with
  the initial data $(\sigma_0(r), R_0)$ satisfying $(\ref{2.1})$. From
  $(\ref{3.7})$, define
\begin{equation}\label{4.1}
  v(r,t):=
   (1-{\gamma\over R(t)}){R(t)\sinh r\over r\sinh R(t)}H(R(t))\qquad\mbox{for }
   0<r\le R(t),\;\;t\ge0.
\end{equation}
  We have the following preliminary lemma.
  
\medskip

{\bf Lemma 4.2} \ \ {\em Let $L>0$ and $M>0$. For some $T>0$, assume that
$$
  |R'(t)|\le L
  \qquad\mbox{for}\;\; 0\le t\le T,
$$
  and 
$$
  |\sigma_0(r)-v(r,0)|\le M \qquad
\mbox{for}\;\; 0\le r\le R_0.
$$
  Then there exists a positive constant $C$ depending only on 
  $\gamma$ such that 
$$
  |\sigma(r,t)-v(r,t)|\le C(Lc+Me^{-{t\over c}})\qquad\mbox{for}\;\;
  0\le r\le R(t), \;\; 0\le t\le T.
$$
}
\medskip 

  {\bf Proof}. \ \  By a direct computation,
$$
  {\partial v\over \partial t}={R(t)\sinh r\over r \sinh R(t)}\Big\{H(R(t))\big[
  {\gamma\over R^2(t)}-R(t)f(R(t))\big]
  +(1-{\gamma\over R(t)})H'(R(t))\Big\}R'(t).
$$
  Since $r f(r)H(r)$ is bounded on $(0,\infty)$ by (\ref{3.4})--(\ref{3.6}), we see that
\begin{equation}\label{4.2}
  |{\partial v\over \partial t}|\le CL
  \qquad\mbox{for}\;\; 0<r\le R(t),\;\;t\ge 0,
\end{equation}
  where $C$ is a constant depending only on $\gamma$.  Let
$$
  \sigma_\pm(r,t)=v(r,t)\pm CLc\pm M e^{-{t\over c}}.
$$
  Then by using (\ref{4.2}) we have
$$
\begin{array}{rl}
\displaystyle  c{\partial \sigma_+\over\partial t}-{1\over r^2}{\partial\over \partial r}(r^2
  {\partial \sigma_+\over \partial r})+\sigma\,&\displaystyle =c{\partial v\over \partial t}
  -\Big[{1\over r^2}{\partial\over \partial r}(r^2{\partial v\over \partial r})-v\Big]
  -Me^{-{t\over c}}+CLc+Me^{-{t\over c}}
  \\ [0.3cm]
  &\displaystyle =c{\partial v\over \partial t}+CLc
  \ge -CLc+CLc\ge0.
\end{array}
$$
  On the other hand, we see that 
$$
  \sigma_+(r,0)=v(r,0)+CLc+M\ge\sigma_0(r)\qquad\mbox{for } 0\le r\le R_0,
$$
  and
$$
 {\partial\sigma_+\over \partial r}(0,t)=0,\qquad
  \sigma_+(R(t),t)=G(t)+CLc+M e^{-{t\over c}}>G(t)\qquad \mbox{for }
  t\ge0,
$$
  where $G(t)=(1-\gamma/R(t))H(R(t))$. Thus by the comparison principle
  of second order parabolic differential equations, we have 
$$
  \sigma_+(r,t)\ge \sigma(r,t)\qquad\mbox{for } 0\le r\le R(t),\;\;0\le t\le T.
$$ 
  Similarly, we have
$$
  \sigma_-(r,t)\le \sigma(r,t)\qquad\mbox{for } 0\le r\le R(t),\;\;0\le t\le T.
$$
  The desired result follows from the above two inequalities. \qquad$\Box$

\medskip
  
  Next, we show that for any given initial data $(\sigma_0(r),R_0)$ satisfying
  (\ref{2.1}), the tumor radius $R(t)$ will be bounded
  for sufficiently small positive $c$.
 
\medskip

{\bf Lemma 4.3} \ \ {\em Let $K$, $\delta>0$ and initial radius $R_0>0$ satisfies one 
  of the following conditions: $(i)$ $\max\{R_0,R_{s2}\}+\delta \le K$;  
  $(ii)$ $R_0+\delta\le K< R_{s1}$. 
   Then there exists a positive constant $c_0$ depending only on $\mu$, 
  $\gamma$, $\tilde\sigma$,  $\delta$, $K$ such that if $0<c\le c_0$, then
$$
  R(t)\le K\qquad\mbox{ for all } \;\;t\ge0.
$$ 
}

{\bf Proof}. \  \ $(i)$ Let $\max\{R_0,R_{s2}\}+\delta \le K$. 
  If the assertion is not true, then there exists $t_0>0$ such that $R(t)<K$ for $0\le t<t_0$
  and  $R(t_0)=K$. It implies that $R'(t_0)\ge0$.

  By (\ref{2.1}) we have $|\sigma_0(r)-v(r,0)|\le \bar\sigma=1$.
  By (\ref{2.3}) we easily get that $|R'(t)|\le L$ for $0\le t\le t_0$, where $L>0$ depends on
  $\mu$, $\tilde\sigma$ and $K$.   It follows from Lemma 4.2 that
$$
  |\sigma(r,t)-v(r,t)|\le C(Lc+e^{-{t\over c}})\qquad\mbox{for}\;\;
  0\le r\le R(t), \;\; 0\le t\le t_0.
$$ 
  Thus by (\ref{1.3}) and (\ref{3.8}) we get
\begin{equation}\label{4.3}
\begin{array}{rl}
  R'(t)\,&\displaystyle \le {\mu\over R^2(t)}\int_0^{R(t)}(v-\tilde\sigma)r^2 dr
  +{\mu\over3}C(Lc+e^{-{t\over c}})R(t)
  \\ [0.5 cm]
  &\displaystyle ={\mu\over3}\Big[\Big(F(R(t))-\tilde\sigma\Big)R(t)+C(Lc+e^{-{t\over c}})
  R(t)\Big].
\end{array}
\end{equation}
  Since $R(t_0)=K>R_{s2}$, by (\ref{3.5}) we see that $F(K)-\tilde\sigma<0$. 
  By taking $t=t_0$ in (\ref{4.3}), we get that for sufficiently small $c>0$, there
  holds $R'(t_0)<0$. It is a contradiction to $R'(t_0)\ge0$, and the assertion holds.
  
  $(ii)$ For the case $R_0+\delta\le K< R_{s1}$, note that we also have $F(K)-\tilde\sigma<0$,
  by a similar argument with a slight modification, we complete the proof.  \qquad$\Box$

\medskip

  Now we study the stability of radially symmetric stationary solution
  $(\sigma_{s2}(r),R_{s2})$. We have the following 
  assertion:

\medskip
  
{\bf Lemma 4.4} \ \ {\em Let $0<\delta<\min\{R_{s2}-R_{s1},1/R_{s2}\}$ 
 and $R_{s1}+\delta<R_0<{1/ \delta}$. For a given $\alpha_0>0$, 
 there exist constants $C$, $b$ and $c_0$ depending on $\mu$, 
 $\gamma$, $\tilde\sigma$, $\delta$, $\alpha_0$ such that if $0<c\le c_0$:
  For any $0<\alpha\le\alpha_0$, if
$$
  |R(t)-R_{s2}|\le \alpha,\qquad
  |R'(t)|\le \alpha,
  \qquad |\sigma(r,t)-\sigma_{s2}(r)|\le\alpha
$$
  hold for all  $0\le r\le R(t)$ and $t\ge0$, then
$$
  |R(t)-R_{s2}|\le C\alpha(c+e^{-bt}),\quad
  |R'(t)|\le C\alpha(c+e^{-bt}),
  \quad |\sigma(r,t)-\sigma_{s2}(r)|\le C\alpha(c+e^{-bt})
$$
  hold for all $0\le r\le R(t)$ and $t\ge T_0$ for some $T_0>0$.
 }

\medskip

{\bf Proof}. \ \ It is easy to verify that there exists a constant 
 $C_0$ depending only on $\gamma$ such that
$$
  |\sigma_0(r)-v(r,0)|\le |\sigma_0(r)-\sigma_{s2}(r)|+|v(r,0)-\sigma_{s2}(r)|
  \le C_0\alpha.
$$
  Then by Lemma 4.2 we have
$$
  |\sigma(r,t)-v(r,t)|\le C_1\alpha(c+e^{-{t\over c}})\qquad\mbox{for}\;\;
  0\le r\le R(t), \;\; t\ge0,
$$ 
  where $C_1$ is also a constant depending only on $\gamma$.
  Similarly as (\ref{4.3}), it follows that
\begin{equation}\label{4.4}
  |R'(t)-{\mu\over3}\big(F(R(t))-\tilde\sigma\big)R(t)|
  \le C_1\alpha\mu R(t)(c+e^{-{t\over c}})\qquad
  \mbox{for}\;\;t\ge0.
\end{equation}
  By using the inequality $e^{-x}\le e^{-1}x^{-1}$ for $x>0$,  we have
\begin{equation}\label{4.5}
  |R'(t)-{\mu\over3}\big(F(R(t))-\tilde\sigma\big)R(t)|
  \le C_2\alpha\mu c R(t)\qquad
  \mbox{for}\;\;t\ge t_0,
\end{equation}
  where $t_0>0$ and $C_2=C_1(1+1/t_0)$. Take $c<1$. By (\ref{4.4}) we also have
\begin{equation}\label{4.6}
\begin{array}{rl}
  |R'(t)| \,&\displaystyle \le {\mu\over3}|F(R(t))-\tilde\sigma|R(t)
  + 2C_1\alpha\mu  R(t)
  \\ [0.5 cm]
  &\displaystyle = {\mu\over3}|F(R(t))-F(R_{s2})|R(t)
  + 2C_1\alpha\mu  R(t)
  \\ [0.3 cm]
  &\le C_3\alpha R(t) 
  \qquad\mbox{for}\;\;0\le t\le t_0,
\end{array}
\end{equation}
  where $C_3=2C_1\mu+\mu\sup_{r>0}\{|F'(r)|/3\}$. Due to $R_{s1}+\delta<R_0<1/\delta$, 
  it gives that
\begin{equation}\label{4.7}
  (R_{s1}+\delta)e^{-C_3\alpha_0 t_0}\le R_0 e^{-C_3\alpha t_0}\le
   R(t_0)\le R_0e^{C_3\alpha t_0}
  \le {1\over\delta}e^{C_3\alpha_0 t_0}.
\end{equation}
  Next, we fix $t_0>0$ such that $(R_{s1}+\delta)e^{-C_3\alpha_0 t_0}>
  R_{s1}+\delta/2$. Consider the following problem:
\begin{equation}\label{4.8}
\left\{
\begin{array}{l}
  \displaystyle {d R^\pm\over dt}={1\over 3}\mu R^\pm(t)
  \Big[F(R^\pm(t))-\tilde\sigma\pm
  3C_2\alpha c \Big]\qquad\mbox{for  } t\ge t_0,
  \\ [0.3cm]
  R^\pm(t_0)=R_0e^{\pm C_3\alpha t_0}.
\end{array}
\right.
\end{equation}
  By (\ref{3.5}) we easily have that there exists $c_0>0$ 
  such that for $0<c\le c_0$,
  the equation
\begin{equation}\label{4.9}
  F(R^\pm)-\tilde\sigma\pm 3C_2\alpha c=0
\end{equation}
  has two positive solutions $R^\pm_{s1}$ and $R^\pm_{s2}$, respectively,
  and satisfy 
$$
 R^+_{s1}<R_{s1}<R^-_{s1}<R_{s1}+\delta/2<R^-_{s2}<R_{s2}<R^+_{s2},
$$
$$
  F'(R^\pm_{s1})>0\qquad\mbox{and}\qquad F'(R^\pm_{s2})<0.
$$
  The constant $c_0$ is dependent
  only on $\mu$, $\gamma$, $\tilde\sigma$, $\delta$ and $\alpha_0$. Besides,
  by mean value theorem, there exists a positive constant  $C_4$ depending only 
  on $\gamma$, $\tilde\sigma$, $\delta$ and $\alpha_0$, such that
\begin{equation}\label{4.10}
  |R^\pm_{s2}-R_{s2}|\le C_4 \alpha c.
\end{equation}
  Hence for the solutions $R^\pm(t)$ of initial value problem (\ref{4.8}), we have
$$
  \lim_{t\to\infty} R^\pm(t)=R^\pm_{s2}.
$$
   Moreover, by a similar argument of (A.25) in [\ref{cui-02}], we can prove that
   there exist constants $C>0$, $b>0$ and $T_0>t_0$ such that
 \begin{equation}\label{4.11}
   |R^\pm(t)-R^\pm_{s2}|\le C_5 \alpha e^{-bt}
   \qquad\mbox{for  } t\ge T_0.
\end{equation}
    By (\ref{4.5}), (\ref{4.7}) and comparison principle of differential equations,
  we have
\begin{equation}\label{4.12}
  R^-(t)\le R(t)\le R^+(t)\qquad \mbox{for } t\ge T_0.
\end{equation}
  From (\ref{4.10})--(\ref{4.12}), we see that there exists a constant $C>0$ such that
  for $t\ge T_0$,
$$
\begin{array}{rl}
  |R(t)-R_{s2}|\,& \le |R^+(t)-R_{s2}|+|R^-(t)-R_{s2}|
  \\ [0.3 cm]
  &\le |R^+(t)-R^+_{s2}|+|R^-(t)-R^-_{s2}|+|R^+_{s2}-R_{s2}|+|R^-_{s2}-R_{s2}|
  \\ [0.3 cm]
  &\le C \alpha (c+e^{-bt}).
\end{array}  
$$
   The other two inequalities follow clearly. \qquad$\Box$
  
\medskip

  With the above preparations, we now state our main result of asymptotic
  behavior.

\medskip

{\bf Theorem 4.5} \ \ {\em Let $0<\tilde\sigma<\theta_*$ and the initial data 
  $(\sigma_0(r),R_0)$ satisfy $(\ref{2.1})$. Suppose that for some small 
  $\varepsilon>0$, the initial radius $R_0$ further satisfies:
 $(i)$ $0< R_0\le R_{s1}-\varepsilon$; or 
 $(ii)$ $R_{s1}+\varepsilon<R_0<1/\varepsilon$.
 Then there exists a positive constant $c_0$ depending only on $\mu$, 
 $\gamma$, $\tilde\sigma$ and $\varepsilon$ such that if $0<c\le c_0$ then
$$
  \lim_{t\to\infty} R(t)=\left\{
  \begin{array}{ll}
  0,\qquad & \mbox{in case } (i),
  \\ [0.2 cm]
  R_{s2}, \qquad & \mbox{in case } (ii).
  \end{array}
  \right.
$$
  }
  
\medskip

{\bf Proof}. \ \ $(i)$ Recall that we have already set $\bar\sigma=1$.
  Let $K=R_{s1}-\varepsilon/2$ and $\delta=\varepsilon/2$. We have
  $R_0+\delta\le K<R_{s1}$. By Lemma 4.3, we see that  there exists $c_1>0$ such that
  for any $0<c\le c_1$, there holds $R(t)\le K$ for $t\ge0$. By (\ref{2.3}) we get 
  $|R'(t)|\le \mu |R(t)|\le \mu K$. 
  Thus by using Lemma 4.2 and similarly as (\ref{4.3}), we have
$$
  R'(t)\le {\mu\over3}\Big[\Big(F(R(t))-\tilde\sigma\Big)R(t)+C(Kc+e^{-{t\over c}})
  R(t)\Big]\qquad\mbox{for  } t\ge0.
$$ 
  By (\ref{3.5}) we see that $F(R(t))-\tilde\sigma\le F(K)-\tilde\sigma:=-2\eta<0$.
  It follows that there exists sufficiently small $c_0>0$ such that for $0<c\le c_0$, 
\begin{equation}\label{4.13}
  R'(t)\le -\mu_0 R(t)\qquad\mbox{for  } t\ge1,
\end{equation}
  where $\mu_0=\mu\eta/3>0$. Hence we have $\lim_{t\to\infty}R(t)=0$ and moreover,
  the convergence is exponentially fast.

  $(ii)$ Set $K=\max\{R_{s2},1/\varepsilon\}+\varepsilon$ and $\delta=\min\{
  \varepsilon, R_{s2}-R_{s1},1/R_{s2}\}$. We see that the conditions of  Lemma 4.3
  $(i)$ and Lemma 4.4 hold. By Lemma 4.3, there exists $c_2>0$ such that for 
  $0<c\le c_2$, we have $R(t)\le K$ for $t\ge0$. Then (\ref{2.3}) implies that 
  $|R'(t)|\le \mu K$. By (\ref{2.2}),  $0\le \sigma(r,t)\le 1$, we have 
  $|\sigma(r,t)-\sigma_{s2}(r)|\le 1$. Let $\alpha_0=(1+\mu)K+1$. 
  By Lemma 4.4, there exist positive constants
  $c_3$, $b$, $C$ and $T_0$ such that for $0<c\le c_3$,  we have
  $|R(t)-R_{s2}|\le C  \alpha_0(c+e^{-bt})\le 2Cc \alpha_0 $ on $[T_0,\infty)$.
  Let $c_0=\min\{c_2, c_3, C/4\}$. Then for $0<c\le c_0$, we have 
  $|R(t)-R_{s2}|\le 2Cc_0\alpha_0\le {1\over 2}\alpha_0$ on $[T_0,\infty)$.
  Hence, by iterating this result over the time intervals $[nT_0,\infty)$
  (as in [\ref{fri-rei-99}]) we get the desired assertion. \qquad$\Box$
  
\medskip

  Note that by Theorem 4.1, in case $\tilde\sigma>\bar\sigma=1$, 
  we have $\lim_{t\to\infty}R(t)=0$ for all $c>0$. For the  case 
  $\theta_*<\tilde\sigma\le 1$, the following result holds:
  
\medskip

{\bf Corollary 4.6} \ \ {\em Let $\theta_*<\tilde\sigma\le 1$. For any given 
  initial data $(\sigma_0(r),R_0)$ satisfying $(\ref{2.1})$, 
  there exists a positive constant $c_0$ depending only on $\mu$, 
  $\gamma$, $\tilde\sigma$ such that for $0<c\le c_0$,  we have 
$$
 \displaystyle  \lim_{t\to\infty} R(t)=0.
$$
  }

{\bf Proof}. \ \ It is similar to the proof of Theorem 4.5 $(i)$. We just need to notice
 that in case $\tilde\sigma>\theta_*$, $F(R)-\tilde\sigma\le \theta_*-\tilde\sigma<0$ 
 for all $R>0$. \qquad $\Box$
 
\medskip
\hskip 1em

\section{Conclusion}
\setcounter{equation}{0}
\hskip 1em

  In this paper, we study a free boundary problem modeling tumor growth 
  with Gibbs-Thomson relation, which is based on the hypothesis that 
  tumor cells on the boundary need consume nutrient for providing energy 
  to maintain the compactness, and the consumption is assumed to be 
  measured by $\gamma/R(t)$, where $\gamma$ is cell-to-cell adhesiveness 
  and $1/R(t)$ represents the mean curvature of tumor boundary with 
  radius $R(t)$ at time $t$.  
  
  An interesting phenomenon induced by 
  Gibbs-Thomson relation is that the model may have two radially
  symmetric stationary solutions, more precisely, 
  there exists $0<\theta_*<1$ depending only on $\gamma$,
  such that for $0<\tilde\sigma<\theta_*$,  problem 
  (\ref{1.1})--(\ref{1.6}) has two radially symmetric stationary solutions.
  It is different from the uniqueness
  of well-studied tumor models by assuming the concentration of nutrient 
  is continuous across the boundary, cf. [\ref{cui-05}, \ref{fri-rei-99}, \ref{xu-16}].  
  
  Our analysis shows that in case
  $0<\tilde\sigma<\theta_*$, for the ratio of the diffusion time scale 
  to the tumor doubling time scale $c$ is sufficiently small, the radially
  symmetric stationary solution with the larger radius is asymptotically stable,
  and the other one with the smaller radius
  is unstable; in case $\tilde\sigma>\theta_*$, the tumor will eventually shrink
  and die for sufficiently small $c$, especially in case $\tilde\sigma>1$,
  the tumor will eventually die for all $c>0$.
     
  Our analysis also implies that the cell-to-cell adhesiveness $\gamma$ may
  have a positive effect on stabilizing the tumor growth.  The larger cell-to-cell 
  adhesiveness, the smaller value $\theta_*$ and radius $R_{s2}$ of the stable
  stationary solution. It indicates that increasing cell-to-cell adhesiveness may
  make the tumor eventually converge to a smaller dormant tumor or 
  die more likely. We hope these analysis may be useful for 
  scientific study and clinical treatment of tumors.

\medskip 
\hskip 1em

  {\bf Acknowledgement.}\ \  This work is partially supported by 
  the PAPD of Jiangsu Higher Education Institutions.

\end{document}